\documentclass[12pt]{article} % Specifies the document class
\usepackage{amsfonts}           % The preamble begins here.
\usepackage{myart2k}
\title{\normalsize{\bf CASCADE CONNECTIONS OF LINEAR SYSTEMS AND
FACTORIZATIONS OF HOLOMORPHIC OPERATOR FUNCTIONS AROUND A MULTIPLE
ZERO IN SEVERAL VARIABLES}}
\author{Dmitriy S. Kalyuzhniy}
\date{}
\newcommand{\nspace}[2]{\ensuremath{{\mathbb{#1}}^{#2}}}
\newcommand{\Hspace}[1]{\ensuremath{\mathcal{#1}}}
\def\ifundefined#1{\expandafter\ifx\csname#1\endcsname\relax}
\providecommand{\comment}[1]{}

\usepackage{amsfonts}
\providecommand{\tthdump}[1]{#1}
\tthdump{}
\newcommand{\Cliff}[2][\comment]{{\ensuremath{%
\mathcal{C}\kern-0.18em\ell(#1,#2)}}}

\ifundefined{qed}
    \DeclareMathSymbol{\qed}{0}{AMSa}{"03}
\fi
\providecommand{\eqref}[1]{\textup{(\ref{#1})}}

\providecommand{\href}[2]{#2}

\begin{document}
\maketitle
\vspace{-1cm}
\begin{abstract}
\noindent
We show that the factorization problem $\theta
(z)=\theta_2(z)\theta_1(z)$ is solvable in the class of Hilbert space
operator-valued functions holomorphic on some neighbourhood of $z=0$ in
$\nspace{C}{N}$ and having  a  zero at $z=0$ (here $\theta (z)$ has a
multiple zero at $z=0$). Such a factorization problem becomes more
complicated if we demand for $\theta (z),\  \theta_1(z)$ and
$\theta_2(z)$ to be Agler--Schur-class functions on the polydisk
$\nspace{D}{N}$ and for the factorization identity to hold in
$\nspace{D}{N}$. In this case we reduce it to the problem on the
existence of a cascade decomposition for certain multiparametric linear
system $\alpha $--a conservative realization of $\theta (z)$, and give
the criterion for its solvability in terms of common invariant subspaces
for the $N$-tuple of main operators  of $\alpha $.
\end{abstract}
\section{Preliminaries}\label{sec:prelim}
In this section we recall the necessary information on multiparametric
linear systems, the Agler--Schur class of operator-valued functions on
the
open unit polydisk $\nspace{D}{N}$, and related realization theorems for
holomorphic operator-valued functions of several complex variables.

In \cite{K2} we have introduced \emph{multiparametric discrete
time-invariant linear dynamical systems} of the form
\begin{equation}\label{eq:n-sys}
\alpha :\left\{\begin{array}{lll}
x(t)&=&\sum _{k=1}^N(A_kx(t-e_k)+B_ku(t-e_k)),\\
y(t)&=&\sum _{k=1}^N(C_kx(t-e_k)+D_ku(t-e_k)),
\end{array}
\right. \quad (t\in\nspace{Z}{N},\ |t|>0)
\end{equation}
where $|t|:=\sum_{k=1}^Nt_k,\ e_k:=(0,\ldots , 0, 1, 0,\ldots ,0)\in
\nspace{Z}{N}$ (here $1$ is on the $k$-th place and zeros are
otherwise); for all $t\in\nspace{Z}{N}$ such that $|t|\geq 0\quad x(t)\
(\in\Hspace{X}),\quad u(t)\ (\in\Hspace{U})$, and for all
$t\in\nspace{Z}{N}$ such that $|t|>0\quad y(t)\
(\in\Hspace{Y})$ are respectively \emph{states},  \emph{input data}
and  \emph{output data}
of $\alpha $, and $\Hspace{X}, \Hspace{U}, \Hspace{Y}$ are
separable Hilbert spaces; for all $k\in\{1,\ldots ,N\} \ A_k, B_k, C_k,
D_k$ are bounded linear operators acting in corresponding pairs of
Hilbert spaces. If one denotes an $N$-tuple of operators $T_k\quad
(k=1,\ldots ,N)$ by $\mathbf{T}:=(T_1,\ldots ,T_N)$ then for
system \eqref{eq:n-sys} one may use the short notation $\alpha =(N;
\mathbf{A}, \mathbf{B}, \mathbf{C}, \mathbf{D};\Hspace{X}, \Hspace{U},
\Hspace{Y})$. If $z=(z_1,\ldots , z_N)\in\nspace{C}{N}$ then we set
$z\mathbf{T}:=\sum _{k=1}^Nz_kT_k$. The
operator-valued function
\begin{equation}\label{eq:n-tf}
\theta _\alpha (z)=z\mathbf{D}+z\mathbf{C}{(I_\Hspace{X}-z\mathbf{A})}^
{-1}z\mathbf{B}
\end{equation}
(here  $I_\Hspace{X}$ is the identity operator on $\Hspace{X}$)
holomorphic on some neighbourhood of $z=0$ in $\nspace{C}{N}$ is called
the \emph{transfer function} of a system $\alpha $ of the form
\eqref{eq:n-sys}. We have the following result \cite{K4}: an arbitrary
function $\theta (z)$ which is holomorphic on some neighbourhood
$\Gamma$
of $z=0$ in $\nspace{C}{N}$ and vanishing at $z=0$, whose values are
from $[\Hspace{U}, \Hspace{Y}]$ (we use the notation
$[\Hspace{U}, \Hspace{Y}]$ for the Banach space of all bounded linear
operators mapping a separable Hilbert space $\Hspace{U}$ into a
separable Hilbert space $\Hspace{Y}$) can be realized as the  transfer
function of some system $\alpha =(N; \mathbf{A}, \mathbf{B}, \mathbf{C},
\mathbf{D};\Hspace{X}, \Hspace{U}, \Hspace{Y})$, i.e. $\theta
(z)=\theta_\alpha (z)$ in some neighbourhood (possibly, smaller than
$\Gamma$) of $z=0$.
We call $\alpha =(N; \mathbf{A}, \mathbf{B}, \mathbf{C},
\mathbf{D};\Hspace{X}, \Hspace{U}, \Hspace{Y})$ a \emph{dissipative}
(resp., \emph{conservative}) \emph{scattering system} if for any
$\zeta\in\nspace{T}{N}$ (the $N$-fold unit torus)
\begin{displaymath}
\zeta\mathbf{G}_\alpha :=\left [\begin{array}{ll}
\zeta\mathbf{A} & \zeta\mathbf{B}\\
\zeta\mathbf{C} & \zeta\mathbf{D}
\end{array}
\right ]\in [\Hspace{X}\oplus\Hspace{U}, \Hspace{X}\oplus\Hspace{Y}]
\end{displaymath}
is a contractive (resp., unitary) operator. Recall that the
\emph{Agler--Schur class} $S_N(\Hspace{U},\Hspace{Y})$ (see \cite{Ag})
consists of all holomorphic functions
$\theta (z)=\sum_{t\in\nspace{Z}{N}_+}\widehat{\theta}_tz^t$
on $\nspace{D}{N}$ with values in $[\Hspace{U},\Hspace{Y}]$
(here $\nspace{Z}{N}_+:=\{ t\in\nspace{Z}{N}:\ t_k\geq 0,\ k=1,\ldots
,N\},\-
z^t:=\prod_{k=1}^Nz_k^{t_k}$) such that
for any separable Hilbert space $\Hspace{H}$,
any $N$-tuple $\mathbf{T}=(T_1,\ldots ,T_N)$ of commuting
contractions on $\Hspace{H}$ and  any positive $r<1$ one has
\begin{equation}\label{eq:agler}
\| \theta (r\mathbf{T})\|  \leq  1
\end{equation}
where
\begin{equation}\label{eq:agler'}
\theta (r\mathbf{T})=\theta (rT_1,\ldots
,rT_N) := \sum_{t\in\nspace{Z}{N}_+}\widehat{\theta _t}\otimes
(r\mathbf{T})^t
\end{equation}
(the convergence of this series  is understood in the
sense of norm in
$[\Hspace{U}\otimes \Hspace{H},\Hspace{Y}\otimes \Hspace{H}]$).
For $N=1$ due to the von Neumann inequality \cite{vN} we have
$S_N(\Hspace{U},\Hspace{Y})=S(\Hspace{U},\Hspace{Y})$, i.e. the Schur
class consisting of all functions holomorphic on the open unit disk
$\mathbb{D}$ with contractive values from $[\Hspace{U},\Hspace{Y}]$. In
\cite{K2} we have proved that the class of transfer functions of
$N$-parametric conservative scattering systems with the input
space $\Hspace{U}$ and the output space $\Hspace{Y}$ coincides with
the subclass $S_N^0(\Hspace{U},\Hspace{Y})$ of
$S_N(\Hspace{U},\Hspace{Y})$ consisting of functions vanishing at $z=0$.
Moreover, the conservative realization $\alpha =(N; \mathbf{A},
\mathbf{B},
\mathbf{C},\mathbf{D};\Hspace{X}, \Hspace{U}, \Hspace{Y})$ of an
arbitrary
function $\theta (z)\in S_N^0(\Hspace{U},\Hspace{Y})$ can be chosen
\emph{closely connected}, i.e. such that $\Hspace{X}_{cc}=\Hspace{X}$
where
\begin{equation}\label{eq:cc-space}
\Hspace{X}_{cc}:=\bigvee_{p,\ k,\
j}p(\mathbf{A},\mathbf{A}^*)(B_k\Hspace{U}+C_j^*\Hspace{Y})
\end{equation}
(here ``$\bigvee_\nu \mathcal{L}_\nu$'' denotes the closure of the
linear span
of subsets $\mathcal{L}_\nu$ in $\Hspace{X}$, $p$ runs over the set of
all monomials in $2N$ non-commuting variables,
$k$ and $j$ run over the set $\{ 1,\ldots ,N\}$). Notice that the close
connectedness of a conservative scattering system $\alpha =(N;
\mathbf{A}, \mathbf{B}, \mathbf{C}, \mathbf{D};\Hspace{X}, \Hspace{U},
\Hspace{Y})$
is equivalent to the condition that the linear pencil
$\zeta\mathbf{A}\ (\zeta\in\nspace{T}{N})$ of
contractive operators is completely non-unitary, i.e. there is no proper
subspace in $\Hspace{X}$ reducing $\zeta\mathbf{A}$
to a unitary operator for each $\zeta\in\nspace{T}{N}$.

\section{Certain factorization problems in several complex variables}
In this section we study factorizations of operator-valued functions
which
are holomorphic on a neighbourhood of some point $z=z_0$ in
$\nspace{C}{N}$ and have a zero at $z=z_0$. Without loss of generality,
one can consider $z_0=0$. For this case we say that such a function
$\theta (z)$, which is not vanishing identically on a neighbourhood of
$z=0$, has a  zero of \emph{multiplicity} $m=m(\theta)$ at $z=0$
if $m$ is the least number of a non-zero term in the expansion of
$\theta
(z)$ in homogeneous polynomials (e.g., see \cite{R2} ).
\begin{prob}\label{prob:1}
Given a function $\theta (z)$ which is holomorphic on the neighbourhood
$\Gamma$ of $z=0$ in $\nspace{C}{N}$, takes values from
$[\Hspace{U},\Hspace{Y}]$, and has a zero of multiplicity $m(\theta )>1$
at $z=0$, find a separable Hilbert space $\Hspace{V}$ and functions
$\theta_1(z),\ \theta_2(z)$ which are holomorphic on some neighbourhoods
of $z=0$, take values from $[\Hspace{U},\Hspace{V}]$ and
$[\Hspace{V},\Hspace{Y}]$ respectively, $\theta_1(0)=0,\ \theta_2(0)=0$,
and
\begin{displaymath}
\theta (z)=\theta_2(z)\theta_1(z)
\end{displaymath}
holds in some neighbourhood (possibly, smaller than $\Gamma$)
 of $z=0$.
\end{prob}
This problem is solvable; in fact, the following more strong statement
is true.
\begin{thm}\label{thm:left}
Let $\theta (z)$ be a function holomorphic on some neighbourhood
$\Gamma$ of $z=0$ in $\nspace{C}{N}$, taking values from
$[\Hspace{U},\Hspace{Y}]$ and having a zero of multiplicity $m=m(\theta
)>0$ at $z=0$. Then there exist separable Hilbert spaces
$\Hspace{Y}^{(0)}=\Hspace{Y},\Hspace{Y}^{(1)},\ldots ,\Hspace{Y}^{(m)}$,
operators $L_k^{(j)}\in [\Hspace{Y}^{(j)},\Hspace{Y}^{(j-1)}]\
(j=1,\ldots ,m;\ k=1,\ldots ,N)$, and a function $\phi (z)$ which is
holomorphic on some neighbourhood of $z=0$, takes values from
$[\Hspace{U},\Hspace{Y}^{(m)}]$, and $\phi (0)\neq 0$, such that
\begin{equation}\label{eq:left}
\theta (z)=z\mathbf{L}^{(1)}\cdots z\mathbf{L}^{(m)}\phi (z)
\end{equation}
holds in some neighbourhood (possibly, smaller than $\Gamma$) of $z=0$.
\end{thm}
\begin{proof}
As it was said in Section~\ref{sec:prelim},  we have $\theta
(z)=\theta_\alpha (z)$ in some neighbourhood of $z=0$ where $\alpha =(N;
\mathbf{A}, \mathbf{B}, \mathbf{C}, \mathbf{D};\Hspace{X}, \Hspace{U},
\Hspace{Y})$ is some $N$-parametric system. If $m=1$ then setting
$\Hspace{Y}^{(1)}:=\Hspace{X}\oplus\Hspace{U},\
L_k^{(1)}:=[C_k^{(1)}\ D_k^{(1)}]\ (k=1,\ldots ,N),\ \phi
(z):=\left[\begin{array}{c}
(I_\Hspace{X}-z\mathbf{A})^{-1}z\mathbf{B}\\
I_\Hspace{U}
\end{array}\right]$, we get from \eqref{eq:n-tf} the equality $\theta
(z)=z\mathbf{L}^{(1)}\phi (z)$, with $\phi (0)\neq 0$, i.e.
\eqref{eq:left} is true.
Let us apply the induction on $m$. Suppose that the statement is true
for $m-1\ (m>1)$. Since $m=m(\theta )>1$ implies $\mathbf{D}=(0,\ldots
,0)$ and \eqref{eq:n-tf} turns into $\theta
(z)=z\mathbf{C}(I_\Hspace{X}-z\mathbf{A})^{-1}z\mathbf{B}$ we can set
$\Hspace{Y}^{(1)}:=\Hspace{X},\ L_k^{(1)}:=C_k\ (k=1,\ldots ,N),\
\widetilde{\theta }(z):=(z\mathbf{A})^{m-2}(I_\Hspace{X}-z\mathbf{A})^{-
1}z\mathbf{B}$ and get $\theta (z)=z\mathbf{L}^{(1)}\widetilde{\theta
}(z)$. Indeed, $z\mathbf{L}^{(1)}\widetilde{\theta
}(z)=z\mathbf{C}(z\mathbf{A})^{m-2}(I_\Hspace{X}-z\mathbf{A})^{-
1}z\mathbf{B}=z\mathbf{C}(I_\Hspace{X}-z\mathbf{A})^{-
1}z\mathbf{B}=\theta (z)$ since
$z\mathbf{C}(z\mathbf{A})^jz\mathbf{B}=0$ identically for $j<m-2$ (the
case $m=2$ is obvious). It is clear that $m(\widetilde{\theta })=m-1$.
By the supposition of the induction there exist separable Hilbert spaces
$\Hspace{Y}^{(2)},\ldots ,\Hspace{Y}^{(m)}$, operators $L_k^{(j)}\in
[\Hspace{Y}^{(j)},\Hspace{Y}^{(j-1)}]\ (j=2,\ldots ,m;\ k=1,\ldots ,N)$,
and a holomorphic function $\phi (z)$ with values from
$[\Hspace{U},\Hspace{Y}^{(m)}]$ such that $\phi (0)\neq 0$ and
$\widetilde{\theta }(z)=z\mathbf{L}^{(2)}\cdots z\mathbf{L}^{(m)}\phi
(z) $ in some neighbourhood of $z=0$. Then \eqref{eq:left} is true, and
the  proof is complete.
\end{proof}
\begin{cor}\label{cor:pol}
If $\theta (z)$ is a homogeneous polynomial of degree $m$ then in the 
statement of Theorem~\ref{thm:left} one can
choose $\Hspace{Y}^{(m)}=\Hspace{U},\ \phi (z)=I_\Hspace{U}$, and
\eqref{eq:left} turns into
\begin{equation}\label{eq:pol}
\theta (z)=z\mathbf{L}^{(1)}\cdots z\mathbf{L}^{(m)}.
\end{equation}
If $m>1$ then,  moreover,
\begin{equation}\label{eq:pol'}
\theta (z)=z\mathbf{C}(z\mathbf{A})^{m-2}z\mathbf{B}.
\end{equation}
\end{cor}
\begin{proof}
The case $m=1$ is trivial. For $m>1$  we obtain \eqref{eq:pol'} from
\eqref{eq:n-tf}  by virtue of the uniqueness of Maclaurin's expansion
for $\theta (z)$.
\end{proof}
Similarly, one can obtain the right-hand analogue of
Theorem~\ref{thm:left}.
\begin{thm}\label{thm:right}
Let $\theta (z)$ be a function holomorphic on some neighbourhood
$\Gamma$ of $z=0$ in $\nspace{C}{N}$, taking values from
$[\Hspace{U},\Hspace{Y}]$ and having a zero of multiplicity $m=m(\theta
)>0$ at $z=0$. Then there exist separable Hilbert spaces
$\Hspace{U}^{(0)}=\Hspace{U},\Hspace{U}^{(1)},\ldots ,\Hspace{U}^{(m)}$,
operators $R_k^{(j)}\in [\Hspace{U}^{(j-1)},\Hspace{U}^{(j)}]\
(j=1,\ldots ,m;\ k=1,\ldots ,N)$, and a function $\psi (z)$ which is
holomorphic on some neighbourhood of $z=0$, takes values from
$[\Hspace{U}^{(m)},\Hspace{Y}]$, and $\psi (0)\neq 0$, such that
\begin{equation}\label{eq:right}
\theta (z)=\psi (z)z\mathbf{R}^{(m)}\cdots z\mathbf{R}^{(1)}
\end{equation}
holds in some neighbourhood (possibly, smaller than $\Gamma$) of $z=0$.
\end{thm}
Let us remark that Theorems~\ref{thm:left} and \ref{thm:right} are
multivariate  generalizations of the theorem on a multiple zero for
functions of one complex variable which are different, even for the
scalar-valued case, from the celebrated Weierstrass Preparation Theorem
(WPT). We have in \eqref{eq:left} and \eqref{eq:right} the products of
linear factors (i.e., the homogeneous polynomials) instead of the
Weierstrass polynomial in WPT which is, in fact, a polynomial in one
distinguished variable and not necessarily a polynomial in other
variables (see, e.g., \cite{R2}). However, these linear factors can be
operator-valued, and, in contrast to WPT, factorizations \eqref{eq:left}
and \eqref{eq:right} are non-unique. Let us remark also that for $N=1$
factorizations
\eqref{eq:pol} and \eqref{eq:pol'} are trivial and reduced
to $\theta (z)=z^mL$ where $L\in [\Hspace{U},\Hspace{Y}]$.
\begin{prob}\label{prob:2}
Given a function $\theta (z)\in S_N^0(\Hspace{U}, \Hspace{Y})$ such that
$m(\theta )>1$, find a separable Hilbert space $\Hspace{V}$ and
functions $\theta_1(z)\in S_N^0(\Hspace{U},\Hspace{V}),\ \theta_2(z)\in
S_N^0(\Hspace{V},\Hspace{Y})$ such that
\begin{equation}\label{eq:agler-fact}
\theta (z)=\theta_2(z)\theta_1(z).\quad (z\in\nspace{D}{N})
\end{equation}
\end{prob}
For the special case of a homogeneous polynomial Problem~\ref{prob:2} is
solvable.
\begin{thm}\label{thm:agler-pol}
If $\theta (z)\in S_N^0(\Hspace{U},\Hspace{Y})$ is a homogeneous
polynomial of degree $m$ then \eqref{eq:pol} holds for
$z\in\nspace{D}{N}$ with linear factors $z\mathbf{L}^{(j)}\in
S_N^0(\Hspace{Y}^{(j)},\Hspace{Y}^{(j-1)})\ (j=1,\ldots ,m)$ (here
$\Hspace{Y}^{(0)}=\Hspace{Y}$ and $\Hspace{Y}^{(m)}=\Hspace{U}$).
\end{thm}
\begin{proof}
As it was said in Section~\ref{sec:prelim} there exists a conservative
scattering system $\alpha =(N; \mathbf{A}, \mathbf{B}, \mathbf{C},
\mathbf{D};\Hspace{X}, \Hspace{U}, \Hspace{Y})$ such that $\theta
(z)=\theta_\alpha (z)$ in $\nspace{D}{N}$. If $m=1$ then $\theta
(z)=\theta_\alpha (z)=z\mathbf{D}$, thus the statement is valid with
$\mathbf{L}^{(1)}=\mathbf{D}$. If $m>1$ then $\mathbf{D}=(0,\ldots ,0)$
and \eqref{eq:pol'} holds. Let us show that $z\mathbf{B}\in
S_N^0(\Hspace{U},\Hspace{X})$. Indeed, for the linear operator-valued
function $z\mathbf{G}_\alpha =\left[\begin{array}{cc}
z\mathbf{A} & z\mathbf{B} \\
z\mathbf{C} & z\mathbf{D}
\end{array}\right]$ corresponding to a conservative scattering system
$\alpha $ we have proved in \cite{K3} that $z\mathbf{G}_\alpha \in
S_N^0(\Hspace{X}\oplus\Hspace{U},\Hspace{X}\oplus\Hspace{Y})$, hence for
any $N$-tuple $\mathbf{T}=(T_1,\ldots ,T_N)$ of commuting contractions
on some separable Hilbert space $\Hspace{H}$ we obtain from
\eqref{eq:agler} and \eqref{eq:agler'}
\begin{displaymath}
\| \sum_{k=1}^N(G_\alpha )_k\otimes T_k\|
_{[(\Hspace{X}\oplus\Hspace{U})\otimes
\Hspace{H},(\Hspace{X}\oplus\Hspace{Y})\otimes \Hspace{H}]}\leq 1
\end{displaymath}
(for a finite sum in \eqref{eq:agler'} one can pass to the limit in
\eqref{eq:agler} as $r\uparrow 1$). Then \begin{displaymath}
\|\sum_{k=1}^NB_k\otimes T_k\| =\| P_\Hspace{X}\otimes
I_\Hspace{H}\left.\left(\sum_{k=1}^N(G_\alpha )_k\otimes
T_k\right)\right|\Hspace{U}\otimes \Hspace{H}\|
\leq\|\sum_{k=1}^N(G_\alpha
)_k\otimes T_k\|\leq 1
\end{displaymath}
(here $P_\Hspace{X}$ is the orthoprojector onto $\Hspace{X}$ in
$\Hspace{X}\oplus\Hspace{Y}$), and by virtue of an arbitrariness of
$\Hspace{H}$ and $\mathbf{T}$ we get $z\mathbf{B}\in
S_N^0(\Hspace{U},\Hspace{X})$. Analogously, $z\mathbf{A}\in
S_N^0(\Hspace{X},\Hspace{X}),\ z\mathbf{C}\in
S_N^0(\Hspace{X},\Hspace{Y})$, and by Corollary~\ref{cor:pol} the
statement of Theorem~\ref{thm:agler-pol} is valid for $m>1$ also.
\end{proof}
For the general case Problem~\ref{prob:2} is still  open. However we
shall show how to reformulate this as the problem on the existence of a
cascade decomposition for a conservative
realization of $\theta (z)$  and give the criterion for its solvability
in
terms of common invariant subspaces for the $N$-tuple
$\mathbf{A}=(A_1,\ldots
,A_N)$ of main operators of such a realization.

\section{Cascade connections of multiparametric linear systems and
factorizations of their transfer fun\-ctions}\label{sec:cascade}
In \cite{K5} we have introduced the notion of cascade connection of
systems of the form \eqref{eq:n-sys} and established some their
properties. Recall that for systems
$\alpha^{(1)} =(N; \mathbf{A}^{(1)},\- \mathbf{B}^{(1)}, \mathbf{C}^{(1)},
\mathbf{D}^{(1)};\Hspace{X}^{(1)}, \Hspace{U}, \Hspace{V})$ and
$\alpha^{(2)} =(N; \mathbf{A}^{(2)}, \mathbf{B}^{(2)},
\mathbf{C}^{(2)},
\mathbf{D}^{(2)}; \Hspace{X}^{(2)}, \Hspace{V}, \Hspace{Y})$ their
\emph{cascade connection} is the system $\alpha =\alpha^{(2)}\alpha
^{(1)}= (N;\- \mathbf{A},\- \mathbf{B}, \mathbf{C},
\mathbf{D};\Hspace{X}=\Hspace{X}^{(2)}\oplus\Hspace{V}\oplus\Hspace{X}^{
(1)}, \Hspace{U}, \Hspace{Y})$ where for any $z\in\nspace{C}{N}$
\begin{eqnarray}
z\mathbf{G}_\alpha =\left[\begin{array}{cc}
z\mathbf{A} & z\mathbf{B} \\
z\mathbf{C} & z\mathbf{D}
\end{array}\right] &:=&  \left [\begin{array}{cccc}
z\mathbf{A}^{(2)} & z\mathbf{B}^{(2)} & 0 & 0 \\
0                 & 0 & z\mathbf{C}^{(1)}  & z\mathbf{D}^{(1)} \\
0                 & 0 & z\mathbf{A}^{(1)}  & z\mathbf{B}^{(1)}
\vspace{1mm}     \\
z\mathbf{C}^{(2)} & z\mathbf{D}^{(2)} & 0             & 0
\end{array}
\right ] \begin{picture}(10,10)
\put(-160,-11){\line(1,0){150}}
\put(-47,-27){\line(0,1){60}}
\end{picture} \label{eq:cascade}  \\
& \in &
[\Hspace{X}^{(2)}\oplus\Hspace{V}\oplus\Hspace{X}^{(1)}\oplus\Hspace{U},
\Hspace{X}^{(2)}\oplus\Hspace{V}\oplus\Hspace{X}^{(1)}\oplus\Hspace{Y}].
\nonumber
\end{eqnarray}
Note that systems \eqref{eq:n-sys} have a unit delay, thus in contrast
to the notion of cascade connection of systems without delay (see, e.g.,
\cite{BC} for the case $N=1$) the state space $\Hspace{X}$ of $\alpha $
contains an additional component--the intermediate space $\Hspace{V}$
(see \cite{K5} for details). If both $\alpha ^{(1)}$ and $\alpha ^{(2)}$
are dissipative (resp., conservative) scattering systems then $\alpha
=\alpha ^{(2)}\alpha ^{(1)}$ is also a dissipative (resp., conservative)
scattering system. If $\theta_{\alpha ^{(1)}}(z)$ and $\theta_{\alpha
^{(2)}}(z)$ are holomorphic on some neighbourhood $\Gamma$ of $z=0$ then
such is $\theta_\alpha (z)$ and
\begin{equation}\label{eq:factor-tf}
\theta_\alpha (z)=\theta_{\alpha ^{(2)}\alpha
^{(1)}}(z)=\theta_{\alpha^{(2)}}(z)\theta_{\alpha^{(1)}}(z).\quad
(z\in\Gamma )
\end{equation}
\begin{thm}\label{thm:clos-con}
If $\alpha =\alpha ^{(2)}\alpha ^{(1)}$ is a closely connected system
then both $\alpha ^{(1)}$ and $\alpha ^{(2)}$ are also closely
connected.
\end{thm}
\begin{proof}
It follows from \eqref{eq:cc-space} that $\Hspace{X}_{cc}$ is the
minimal subspace  in $\Hspace{X}$ containing $B_k\Hspace{U},\
C_k^*\Hspace{Y}$ and reducing $A_k$ for all $k\in\{ 1,\ldots ,N\}$. By
the assumption, $\Hspace{X}_{cc}=\Hspace{X}$. If
$\Hspace{X}_{cc}^{(1)}\neq\Hspace{X}^{(1)}$ then by \eqref{eq:cascade}
the subspace
$\Hspace{X}^{(2)}\oplus\Hspace{V}\oplus\Hspace{X}_{cc}^{(1)}\
(\neq\Hspace{X})$ contains  $B_k\Hspace{U},\ C_k^*\Hspace{Y}$ and
reduces
$A_k$ for all $k\in\{ 1,\ldots ,N\}$, that contradicts to the
assumption.
Hence, $\alpha ^{(1)}$ is closely connected. Analogously, $\alpha
^{(2)}$
is closely connected.
\end{proof}
Note that the same is true for systems
without delay for $N=1$ (see \cite{BC}). The converse statement is false
even
for $N=1$.
\begin{example}
Let $l^2=\bigoplus_{n=-\infty}^{+\infty
}\mathbb{C},\ l^2_+=\bigoplus_{n=0}^{+\infty }\mathbb{C},\
l^2_-=\bigoplus_{n=-\infty}^{-1}\mathbb{C}$ be Hilbert spaces of
sequences.  Clearly, $l^2=l^2_+\oplus l^2_-$. Let $U:l^2\to l^2$ be the
two-sided shift operator:
\begin{displaymath}
U:\mbox{col}(\ldots, c_{-1}, \fbox{\ensuremath{c_0}}, c_1,\ldots
)\mapsto\mbox{col}(\ldots, c_{-2}, \fbox{\ensuremath{c_{-1}}},
c_0,\ldots ).
\end{displaymath}
Obviously, $U$ is unitary. Define the systems $\alpha^{(j)} =(1;
A^{(j)},
B^{(j)}, C^{(j)}, D^{(j)}; l^2_+, l^2_-,\- l^2_-),\ j=1,2$, where
\begin{displaymath}
\left[\begin{array}{cc}
A^{(1)} & B^{(1)} \\
C^{(1)} & D^{(1)}
\end{array}\right]=U^{-1}=U^*=\left[\begin{array}{cc}
A^{(2)} & B^{(2)} \\
C^{(2)} & D^{(2)}
\end{array}\right]^*\in [l^2_+\oplus l^2_-,l^2_+\oplus l^2_-].
\end{displaymath}
In particular,
\begin{eqnarray*}
A^{(1)}:l^2_+\to l^2_+, & \mbox{col}(c_0,c_1,\ldots
)\mapsto\mbox{col}(c_1,c_2,\ldots ), \\
C^{(1)}:l^2_+\to l^2_-, & \mbox{col}(c_0,c_1,\ldots
)\mapsto\mbox{col}(\ldots , 0, 0, c_0), \\
A^{(2)}:l^2_+\to l^2_+, & \mbox{col}(c_0,c_1,\ldots
)\mapsto\mbox{col}(0, c_0, c_1, \ldots ), \\
B^{(2)}:l^2_-\to l^2_+, & \mbox{col}(\ldots ,c_{-2},c_{-
1})\mapsto\mbox{col}(c_{-1}, 0, 0,\ldots ).
\end{eqnarray*}
It is clear that both $\alpha ^{(1)}$ and $\alpha ^{(2)}$ are
conservative scattering systems. Moreover, they are closely connected
since $A^{(1)}$ and $A^{(2)}$ are respectively the backward shift and
the
forward shift operators on $l^2_+$ which are completely non-unitary. Let
$\alpha =\alpha ^{(2)}\alpha ^{(1)}$. Then by \eqref{eq:cascade} the
main operator of $\alpha $ is
\begin{displaymath}
A=\left[\begin{array}{ccc}
A^{(2)} & B^{(2)} & 0 \\
0       & 0       & C^{(1)} \\
0       & 0       & A^{(1)}
\end{array}\right]:l^2_+\oplus l^2_-\oplus l^2_+\to l^2_+\oplus l^2_-
\oplus l^2_+.
\end{displaymath}
Let $\Hspace{K}:=\{ \mbox{col}(\ldots , 0, 0, c_{-1}),\ c_{-
1}\in\mathbb{C}\}\subset l^2_-$, and
$\Hspace{X}_\Hspace{K}:=l^2_+\oplus\Hspace{K}\oplus l^2_+\subset
l^2_+\oplus l^2_-\oplus l^2_+$. Then $A$ is acting on elements of
$\Hspace{X}_\Hspace{K}$ as follows:
\begin{displaymath}
A:\left[\begin{array}{l}
\mbox{col}(c_0^{(2)}, c_1^{(2)}, c_2^{(2)} \ldots ) \\
\mbox{col}(\ldots , 0, 0, c_{-1}) \\
\mbox{col}(c_0^{(1)}, c_1^{(1)}, c_2^{(1)} \ldots )
\end{array}\right]\mapsto
\left[\begin{array}{l}
\mbox{col}(c_{-1}, c_0^{(2)}, c_1^{(2)}, \ldots ) \\
\mbox{col}(\ldots , 0, 0, c_0^{(1)}) \\
\mbox{col}(c_1^{(1)}, c_2^{(1)}, c_3^{(1)}, \ldots )
\end{array}\right].
\end{displaymath}
It is clear now that $\Hspace{X}_\Hspace{K}$ is invariant subspace for
$A$, and $A|\Hspace{X}_\Hspace{K}$ is unitary, thus operator $A$
has a unitary part. Therefore the conservative scattering system $\alpha
$ is not closely connected.
\end{example}
Note that the analogous (however, more complicated)
example was constructed in \cite{Br} for the
case of one-parametric conservative scattering systems without delay
(in the language of unitary colligations).
\begin{thm}\label{thm:casc-decomp}
1. Let $\alpha =(N; \mathbf{A}, \mathbf{B}, \mathbf{C},
\mathbf{D};\Hspace{X}, \Hspace{U}, \Hspace{Y})$ be a conservative
scattering system, and $\Hspace{X}^{(2)}$ be a subspace of $\Hspace{X}$
satisfying the following conditions:

(i) $\Hspace{X}^{(2)}$ is invariant for all $A_k,\ k=1,\ldots ,N$
(equivalently, $\forall\zeta\in\nspace{T}{N}\  \zeta\mathbf{G}_\alpha
\Hspace{X}^{(2)}\subset\Hspace{X}^{(2)}\oplus\Hspace{Y}$);

(ii) the subspaces $\Hspace{V}_\zeta :=(\zeta\mathbf{G}_\alpha
)^*(\Hspace{X}^{(2)}\oplus\Hspace{Y})\ominus\Hspace{X}^{(2)}$ coincide
for all $\zeta\in\nspace{T}{N}\  (\forall\zeta\in\nspace{T}{N}\
\Hspace{V}_\zeta =:\Hspace{V}$).

Define $\Hspace{X}^{(1)}:=\Hspace{X}\ominus
(\Hspace{X}^{(2)}\oplus\Hspace{V})$, for all $k\in\{ 1,\ldots ,N\}$
\begin{eqnarray}
(G_{\alpha ^{(1)}})_k:=P_{\Hspace{X}^{(1)}\oplus\Hspace{V}}(G_\alpha
)_k|\Hspace{X}^{(1)}\oplus\Hspace{U}, \nonumber \\
(G_{\alpha ^{(2)}})_k:=P_{\Hspace{X}^{(2)}\oplus\Hspace{Y}}(G_\alpha
)_k|\Hspace{X}^{(2)}\oplus\Hspace{V}, \label{eq:decomp}
\end{eqnarray}
i.e.,
\begin{eqnarray*}
\left[\begin{array}{cc}
A_k^{(1)} & B_k^{(1)} \\
C_k^{(1)} & D_k^{(1)}
\end{array}\right] &:=& \left[\begin{array}{cc}
P_{\Hspace{X}^{(1)}\oplus\Hspace{V}} & 0 \\
0                                    & 0
\end{array}\right]\left.\left[\begin{array}{cc}
A_k & B_k \\
C_k & D_k
\end{array}\right]\right|\Hspace{X}^{(1)}\oplus\Hspace{U}, \\
\left[\begin{array}{cc}
A_k^{(2)} & B_k^{(2)} \\
C_k^{(2)} & D_k^{(2)}
\end{array}\right] &:=& \left[\begin{array}{cc}
P_{\Hspace{X}^{(2)}} & 0 \\
0                    & I_\Hspace{Y}
\end{array}\right]\left.\left[\begin{array}{cc}
A_k & B_k \\
C_k & D_k
\end{array}\right]\right|(\Hspace{X}^{(2)}\oplus\Hspace{V})\oplus\{ 0\}.
\end{eqnarray*}
Then $\alpha =\alpha ^{(2)}\alpha ^{(1)}$ where $\alpha^{(1)} =(N;
\mathbf{A}^{(1)}, \mathbf{B}^{(1)}, \mathbf{C}^{(1)},
\mathbf{D}^{(1)};\Hspace{X}^{(1)}, \Hspace{U}, \Hspace{V}),\
\alpha^{(2)} =(N;\- \mathbf{A}^{(2)}, \mathbf{B}^{(2)}, \mathbf{C}^{(2)},
\mathbf{D}^{(2)};\Hspace{X}^{(2)}, \Hspace{V}, \Hspace{Y});\ \alpha
^{(1)}$ and $\alpha ^{(2)}$ are conservative scattering systems.

2. Any cascade connection of conservative scattering systems of the form
\eqref{eq:n-sys} arises in
this way.
\end{thm}
\begin{proof}
For any $\zeta\in\nspace{T}{N}\ \zeta\mathbf{G}_\alpha$ is a unitary
operator, and from \eqref{eq:decomp} and the definition of $\Hspace{V}$
we have $\zeta\mathbf{G}_\alpha
(\Hspace{X}^{(2)}\oplus\Hspace{V})=\Hspace{X}^{(2)}\oplus\Hspace{Y}=
\zeta\mathbf{G}_{\alpha ^{(2)}}(\Hspace{X}^{(2)}\oplus\Hspace{V})$, and
$\zeta\mathbf{G}_{\alpha ^{(2)}}=\zeta\mathbf{G}_\alpha
|\Hspace{X}^{(2)}\oplus\Hspace{V}$ is an isometry, thus
$\zeta\mathbf{G}_{\alpha ^{(2)}}$ is unitary. Hence, $\alpha ^{(2)}$ is
a conservative scattering system. For any $\zeta\in\nspace{T}{N}$
\begin{eqnarray*}
\zeta\mathbf{G}_\alpha (\Hspace{X}^{(1)}\oplus\Hspace{U}) &=&
\zeta\mathbf{G}_\alpha
(\Hspace{X}\oplus\Hspace{U})\ominus\zeta\mathbf{G}_\alpha
(\Hspace{X}^{(2)}\oplus\Hspace{V})=(\Hspace{X}\oplus\Hspace{Y})\ominus
(\Hspace{X}^{(2)}\oplus\Hspace{Y}) \\
&=&
\Hspace{X}\ominus\Hspace{X}^{(2)}=\Hspace{X}^{(1)}\oplus\Hspace{V}=\zeta
\mathbf{G}_{\alpha ^{(1)}}(\Hspace{X}^{(1)}\oplus\Hspace{U}),
\end{eqnarray*}
and $\zeta\mathbf{G}_{\alpha ^{(1)}}=\zeta\mathbf{G}_\alpha
|\Hspace{X}^{(1)}\oplus\Hspace{U}$ is an isometry, thus
$\zeta\mathbf{G}_{\alpha ^{(1)}}$ is unitary. Hence, $\alpha ^{(1)}$ is
a
conservative scattering system. It is easy to see now that for any
$\zeta\in\nspace{T}{N}\ \zeta\mathbf{G}_\alpha $ has a form
\eqref{eq:cascade}, i.e. $\alpha =\alpha ^{(2)}\alpha ^{(1)}$. The
second assertion of this theorem follows directly from
the definition of cascade connection.
\end{proof}
Note that Theorem~\ref{thm:casc-decomp} is an analogue of the
well-known result for one-parametric systems without delay (Theorem~2.6
in \cite{Br}; see also
Theorem~6.1 in \cite{BC}).
\begin{thm}\label{thm:fact}
For $\theta (z)\in S_N^0(\Hspace{U},\Hspace{Y})$ such that $m(\theta
)>1$ Problem~\ref{prob:2} is solvable if and only if there exists  a
closely connected conservative scattering system $\alpha_{cc}=(N;
\mathbf{A}_{cc}, \mathbf{B}_{cc}, \mathbf{C}_{cc},
\mathbf{D}_{cc};\Hspace{X}_{cc}, \Hspace{U}, \Hspace{Y})$ such that
$\theta (z)=\theta_{\alpha_{cc}} (z)$ for $z\in\nspace{D}{N}$ and there
exists a subspace $\Hspace{X}^{(2)}_{cc}$ in $\Hspace{X}_{cc}$
satisfying
conditions (i), (ii) in Theorem~\ref{thm:casc-decomp} applied for the
system $\alpha _{cc}$ in the place of $\alpha $.
\end{thm}
\begin{proof}
The part ``if'' is clear since in this case by
Theorem~\ref{thm:casc-decomp} there are conservative scattering systems
$\alpha ^{(1)}$ and $\alpha ^{(2)}$ such that $\alpha _{cc}=\alpha
^{(2)}\alpha ^{(1)}$, and by \eqref{eq:factor-tf} we have $\theta
(z)=\theta_{\alpha _{cc}}(z)=\theta_{\alpha ^{(2)}}(z)\theta_{\alpha
^{(1)}}(z)$ with functions $\theta_{\alpha ^{(1)}}(z)$ and
$\theta_{\alpha ^{(2)}}(z)$ belonging to the corresponding
classes $S_N^0(\cdot ,\cdot )$.

For the proof of the part ``only if'' let us assume that
\eqref{eq:agler-fact} holds with $\theta_1(z)\in
S_N^0(\Hspace{U},\Hspace{V}),\ \theta_2(z)\in
S_N^0(\Hspace{V},\Hspace{Y})$. Let $\alpha^{(1)} =(N; \mathbf{A}^{(1)},
\mathbf{B}^{(1)}, \mathbf{C}^{(1)}, \mathbf{D}^{(1)};\-
\Hspace{X}^{(1)},
\Hspace{U}, \Hspace{V}),\ \alpha^{(2)} =(N; \mathbf{A}^{(2)},
\mathbf{B}^{(2)}, \mathbf{C}^{(2)}, \mathbf{D}^{(2)};\Hspace{X}^{(2)},
\Hspace{V}, \Hspace{Y})$ be some conservative realizations, respectively
(which exist by \cite{K2}), i.e. $\theta_k(z)=\theta_{\alpha
^{(k)}}(z),\ k=1,2$. Then the conservative scattering system $\alpha
=(N; \mathbf{A}, \mathbf{B}, \mathbf{C}, \mathbf{D};\Hspace{X},
\Hspace{U}, \Hspace{Y}):=\alpha ^{(2)}\alpha^{(1)}$ has the transfer
function $\theta_\alpha (z)=\theta_{\alpha ^{(2)}}(z)\theta_{\alpha
^{(1)}}(z)$, and the subspace $\Hspace{X}^{(2)}$ in $\Hspace{X}$ satisfy
conditions (i) and (ii) of Theorem~\ref{thm:casc-decomp}. Define
the subspace $\Hspace{X}_{cc}$ in $\Hspace{X}$ by \eqref{eq:cc-space},
and operators $(A_{cc})_k:=P_{\Hspace{X}_{cc}}A_k|\Hspace{X}_{cc},\
(B_{cc})_k:=P_{\Hspace{X}_{cc}}B_k,\ (C_{cc})_k:=C_k|\Hspace{X}_{cc},\
(D_{cc})_k:=D_k=0,\ k=1,\ldots ,N$. Then, by Theorem~3.3 of \cite{K2},
$\Hspace{X}_{cc}$ is a reducing subspace in $\Hspace{X}$ for all $A_k,\
k=1,\ldots ,N$, and $\alpha_{cc} =(N; \mathbf{A}_{cc}, \mathbf{B}_{cc},
\mathbf{C}_{cc}, \mathbf{D}_{cc};\Hspace{X}_{cc}, \Hspace{U},
\Hspace{Y})$ is a closely connected conservative realization of
$\theta(z)=\theta_\alpha (z)$. Define
$\Hspace{X}_{cc}^{(2)}:=\overline{P_{\Hspace{X}_{cc}}\Hspace{X}^{(2)}}$.
Then $\Hspace{X}_{cc}^{(2)}$ is an invariant subspace in
$\Hspace{X}$ for all $A_k,\ k=1,\ldots ,N$. Since
$(A_{cc})_k=A_k|\Hspace{X}_{cc}$ we obtain that $\Hspace{X}_{cc}^{(2)}$
is an invariant subspace in $\Hspace{X}_{cc}$ for all $(A_{cc})_k,\
k=1,\ldots ,N$, i.e. condition (i) in Theorem~\ref{thm:casc-decomp} is
satisfied for $\Hspace{X}_{cc}^{(2)}$. For all $\zeta\in\nspace{T}{N}$
we have the spaces
\begin{eqnarray*}
(\Hspace{V}_{cc})_\zeta &:=& (\zeta\mathbf{G}_{\alpha
_{cc}})^*(\Hspace{X}_{cc}^{(2)}\oplus\Hspace{Y})\ominus\Hspace{X}_{cc}^{
(2)} \\
&=& (\zeta\mathbf{G}_{\alpha
_{cc}})^*(\overline{P_{\Hspace{X}_{cc}}\Hspace{X}^{(2)}}\oplus\Hspace{Y}
)\ominus\overline{P_{\Hspace{X}_{cc}}\Hspace{X}^{(2)}} \\
 &=& (\zeta\mathbf{G}_\alpha )
^*\overline{P_{\Hspace{X}_{cc}\oplus\Hspace{Y}}(\Hspace{X}^{(2)}\oplus
\Hspace{Y})}\ominus\overline{P_{\Hspace{X}_{cc}}\Hspace{X}^{(2)}} \\
&=& \overline{P_{(\zeta\mathbf{G}_\alpha
)^*(\Hspace{X}_{cc}\oplus\Hspace{Y})}(\zeta\mathbf{G}_
\alpha
)^*(\Hspace{X}^{(2)}\oplus\Hspace{Y})}\ominus\overline{P_{\Hspace{X}_{cc
}}\Hspace{X}^{(2)}} \\
&=& \overline{P_{\Hspace{X}_{cc}\oplus\Hspace{U}}(\Hspace{X}^{(2)}\oplus
\Hspace{V})}\ominus\overline{P_{\Hspace{X}_{cc}}\Hspace{X}^{(2)}} \\
&=& \overline{P_{\Hspace{X}_{cc}}(\Hspace{X}^{(2)}\oplus
\Hspace{V})}\ominus\overline{P_{\Hspace{X}_{cc}}\Hspace{X}^{(2)}}\
(=:\Hspace{V}_{cc})
\end{eqnarray*}
coinciding, and condition (ii) in Theorem~\ref{thm:casc-decomp} is also
satisfied for $\Hspace{X}_{cc}^{(2)}$. The proof is complete.
\end{proof}
Let us give some remarks. {\bf 1}. The subspace $\Hspace{X}_{cc}^{(2)}$
in Theorem~\ref{thm:fact} corresponds to the factorization of $\theta (z)$
which, in general, not necessarily coincides with the original
factorization \eqref{eq:agler-fact}. {\bf 2}. The case when
$\Hspace{X}_{cc}^{(2)}=\{ 0\}$ is also non-trivial (!) since for that
$\Hspace{V}_{cc}=\overline{P_{\Hspace{X}_{cc}}\Hspace{V}}=(\zeta\mathbf{
G}_
{\alpha
_{cc}})^*\Hspace{Y}=(\zeta\mathbf{G}_\alpha )^*\Hspace{Y}\neq \{ 0\}$,
and the factor from the left in the corresponding factorization of
$\theta (z)$ is a linear homogeneous operator-valued function. {\bf 3}.
It is easy to show that one can define $\Hspace{X}_{cc}^{(2)}$ in
another way, say,
$\Hspace{X}_{cc}^{(2)}:=\Hspace{X}_{cc}\cap\Hspace{X}^{(2)}$, and this
subspace also satisfies conditions (i) and (ii) in
Theorem~\ref{thm:casc-decomp}.

\bibliographystyle{amsplain}
\bibliography{msys}
\ \\
Department of Higher Mathematics \\
Odessa State Academy of Civil Engineering and Architecture \\
Didrihson str. 4, Odessa, 270029, Ukraine
\end{document}